\documentclass[11pt,twoside]{article} 
\usepackage[cp1251]{inputenc}
\usepackage[russian]{babel}
\usepackage{amsmath,amsfonts,amssymb}
\usepackage{latexsym,hhline,graphics,epsfig}

\begin{document}

\vspace*{-10mm}

\vskip 3mm

\pagestyle{myheadings} 
\markboth {\underline{\mbox{\!\!}\hspace{6,2cm} \small О.К. Бахтін,
Я.В. Заболотний}}{\underline{\hspace{0.8cm}\small Застосування
розділяючого перетворення... \hfill}}

\begin{flushleft}

\noindent {\small УДК\ \ 517.5}

\bigskip

{\large\bf О.К. Бахтін, Я.В. Заболотний}

\vskip 1mm

\bigskip

{\LARGE\bf Оцінки добутку внутрішніх радіусів  \\[2mm] чотирьох взаємно неперетинних областей}

\renewcommand{\thefootnote}{}
\footnote{\hfill  \copyright \ \ {\bf О.К. Бахтін, Я.В. Заболотний,
2015}}

\vskip 2mm

\end{flushleft}

\vskip 2mm

\begin{flushright}
\parbox{10.5cm}{\footnotesize
В даній роботі розглядається відома гіпотеза В.М. Дубиніна про
неперетинні області на комплексній площині і знайдено її частковий
розв'язок.}
\end{flushright}

\medskip

\vskip 3mm

\thispagestyle{empty} В геометричній теорії функцій комплексної
змінної значне місце займають екстремальні задачі на класах
областей, що не перетинаються. В даному напрямку відзначимо роботи
[1-13]. Було отримано багато вагомих результатів, але в той же час
значна кількість задач не розв'язана і досі. Одній з таких задач і
присвячена дана робота. Дана задача була сформульована в роботі [7,
стор.68].

 Нехай $\mathbb{N}$, $\mathbb{C}$ -- множини
натуральних і комплексних чисел відповідно, $B_{j}$ - область в
$\overline{\mathbb{C}}$, $r(B_{j},a_{j})$ - внутрішній радіус
області $B_{j}$ в точці $a_{j}$.

\textbf{Задача 1.} Довести, що максимум функціонала

$$I_{\gamma}=r^{\gamma}(B_{0},a_{0})\cdot\prod\limits_{k=1}^n
r(B_{k},a_{k})$$

де $B_{0}$, $B_{1}$, $B_{2}$,...,$B_{n}$, $(n\geq2)$ - попарно
неперетині області в $\overline{\mathbb{C}}$, $a_{0}=0$,
$|a_{k}|=1$, $k=\overline{1,n}$, $r(B_{j},a_{j})$ - внутрішній
радіус області $B_{j}$ в точці $a_{j}$ $(a_{j}\in B_{j})$,
$j=\overline{0,n}$ і $\gamma\leq n$ (див, напр. [2-12]), досягається
для деякої конфігурації областей, які мають $n$ - кратну симетрію. В
даній роботі встановлено наступний результат:

\textbf{Теорема 1.}  При $n=3$ і $\gamma\in(0;1.7]$ максимум
функціонала $I_{\gamma}$ досягається на системі областей $D_{0}$,
$D_{1}$, $D_{2}$ , $D_{3}$ і точках $a_{0}$, $a_{1}$, $a_{2}$,
$a_{3}$ де $D_{k}$, $a_{k}$ - відповідно кругові області і полюси
квадратичного диференціала

\begin{equation}
Q(w)dw^{2} =-
\frac{(9-\gamma)w^{3}+\gamma}{w^{2}(w^{3}-1)^{2}}dw^{2}.
\end{equation}

\vskip 2mm

\textbf{Доведення теореми 1.}

Зауважимо, що в роботі [15] було встановлено даний результат для
$n=3$ і $\gamma\in(0;1.5]$. Тому нам залишилося лише розглянути
випадок $\gamma\in(1.5;1.7]$.

Встановимо спочатку, що дане твердження правильне для $\gamma=1.7$.
Метод доведення спирається на застосування, аналогічно до теореми
5.2.3 роботи [9], методу розділяючого перетворення областей, який
детально розроблений в роботі [7]. Для даного випадку розділяюче
перетворення було проведемо так, як і роботі [15].

Далі, як і в згаданій вище теоремі 5.2.3 [9], отримаємо нерівність:

\begin{equation}
I_{\gamma}\leq \left[\prod\limits_{k=1}^3
\alpha_{k}r^{\gamma\cdot\alpha_{k}^{2}}(D_{0},0)\cdot
r(D_{1},i)\cdot r(D_{2},-i)\right]^{\frac{1}{2}}.
\end{equation}

 де $D_{k}$ - згадані вище кругові області квадратичного диференціала (1). Дана
нерівність правильна для $\gamma\leq1$ на основі результатів роботи
[8]. При $\gamma>1$ її застосування, взагалі кажучи, некоректне. ЇЇ
застосування можливе у випадку $\alpha_{k}\sqrt{\gamma}\leq2$.
Знайдемо умови правильності цієї нерівності для $\gamma=1.7$.

Нехай для конкретності $r(B_{0},0)=p$. Тоді

$$I_{\gamma}=r^{\gamma}(B_{0},a_{0})\cdot\prod\limits_{k=1}^n
r(B_{k},a_{k})=p^{\gamma-1}\cdot\prod\limits_{k=0}^n
r(B_{k},a_{k})\leq$$

\begin{equation}
\leq
p^{\gamma-1}\cdot\frac{9}{4^{\frac{8}{3}}}(|a_{1}-a_{2}|\cdot|a_{1}-a_{3}|\cdot|a_{2}-a_{3}|\cdot|a_{1}-a_{0}|\cdot|a_{2}-a_{0}|\cdot|a_{3}-a_{0}|)^{\frac{2}{3}}
\end{equation}

Остання нерівність правильна на основі теореми Кузьміної [14].

Доведемо, що області, які можуть бути екстремальними, задовольняють
умову $\alpha_{0}\leq\frac{2}{\sqrt{\gamma}}$. Обчислимо значення
функціонала

\begin{center}
$I_{\gamma}^{0}=r^{\gamma}(D_{0},a_{0})\cdot\prod\limits_{k=1}^3
r(D_{k},a_{k})$
\end{center}

Згідно із вище згаданою теоремою 5.2.3 роботи [9]

\begin{center}
$I_{\gamma}^{0}=\left(\frac{4}{n}\right)^{n}\cdot\frac{\left(\frac{4\gamma}{n^{2}}\right)^{\frac{\gamma}{n}}}{\left(1-\frac{\gamma}{n^{2}}\right)^{n+\frac{\gamma}{n}}}\cdot
\left(\frac{1-\frac{\sqrt{\gamma}}{n}}{1+\frac{\sqrt{\gamma}}{n}}\right)^{2\sqrt{\gamma}}$
\end{center}

Підставивши $\gamma=1.7$ і $n=3$, отримаємо, що
$I_{1.7}^{0}\approx0.3763$

Далі нам буде потрібна наступна

\textbf{Лема 1.} Нехай $B_{0}$, $B_{1}$, $B_{2}$,...,$B_{n}$,
$(n\geq2)$ - попарно неперетині області в $\overline{\mathbb{C}}$,
$a_{0}=0$, $|a_{k}|=1$, $k=\overline{1,n}$, $(a_{j}\in B_{j})$,
$j=\overline{0,n}$ і $q>0$, $q\in\mathbb{R}$, $r(B_{j},a_{j})$ -
внутрішній радіус області $B_{j}$ в точці $a_{j}$ $(a_{j}\in B_{j})$
і $\gamma\leq n$. Тоді при умові, що $r(B_{0},a_{0})\geq
q^{\frac{1}{\gamma-n}}$ виконується нерівність:

\begin{center}
$r^{\gamma}(B_{0},a_{0})\cdot\prod\limits_{k=1}^n r(B_{k},a_{k})\leq
q$
\end{center}

\textbf{Доведення.} Нехай $r(B_{0},a_{0})=p\geq
q^{\frac{1}{\gamma-n}}$ Застосуємо теорему Лаврентьєва [1] для
областей $B_{0}$ і $B_{1}$, отримаємо, що $r(B_{0},0)\cdot
r(B_{1},a_{1})\leq \mid a_{1}\mid=1$. Оскільки $r(B_{0},0)=p$, то
$r(B_{1},a_{1})\leq\frac{1}{p}$. Аналогічно
$r(B_{k},a_{k})\leq\frac{1}{p}$ для $k=\overline{1,n}$. Тоді

\begin{center}
$r^{\gamma}(B_{0},a_{0})\cdot\prod\limits_{k=1}^n r(B_{k},a_{k})\leq
p^{\gamma}\cdot\frac{1}{p^{n}}=p^{\gamma-n}\leq
(q^{\frac{1}{\gamma-n}})^{\gamma-n}=q$
\end{center}

Лему доведено.

Далі, при $\gamma=1.7$ і $n=3$ застосуємо тільки що доведену Лему 1,
взявши $q=I_{1.7}^{0}$. Таким чином ми отримаємо, що для
$r(B_{0},a_{0})\geq (I_{1.7}^{0})^{\frac{1}{\gamma-n}}\approx2.1208$
$r^{\gamma}(B_{0},a_{0})\cdot\prod\limits_{k=1}^n r(B_{k},a_{k})\leq
I_{1.7}^{0}$, тому конфігурації областей при таких значеннях
$r(B_{0},a_{0})$ не можуть бути екстремальними.

Нехай тепер $p\leq p_{0}=(I_{1.7}^{0})^{\frac{1}{\gamma-n}}$. Тоді,
за (3),
$$I_{\gamma}\leq p^{\gamma-1}\frac{9}{4^{\frac{8}{3}}}\left(|a_{1}-a_{2}|\cdot|a_{1}-a_{3}|\cdot|a_{2}-a_{3}|\right)^{\frac{2}{3}}.$$
Далі, нехай $\alpha_{0}\geq\frac{2}{\sqrt{\gamma}}$. Візьмемо для
конкретності $\alpha_{1}=\alpha_{0}$. Тоді

$$|a_{1}-a_{2}|=2\sin\frac{\alpha_{1}\cdot\pi}{2}\leq2\sin\frac{2-\frac{2}{\sqrt{\gamma}}\cdot\pi}{2}=2\sin\left(1-\frac{1}{\sqrt{\gamma}}\right)\cdot\pi=2\cdot\sin\frac{\pi}{\sqrt{\gamma}}$$

Далі, оскільки
$\alpha_{2}+\alpha_{3}=2-\alpha_{1}\leq\frac{2}{\sqrt{\gamma}}$, за
нерівністю Коші максимальне значення добутку
$|a_{1}-a_{3}|\cdot|a_{2}-a_{3}|$ ми отримаємо в тому випадку, коли
$|a_{1}-a_{3}|=|a_{2}-a_{3}|$, тобто при
$\alpha_{2}=\alpha_{3}=1-\frac{1}{\sqrt{\gamma}}$. Звідси
$|a_{1}-a_{3}|=|a_{2}-a_{3}|=2\sin\left(1-\frac{1}{\sqrt{\gamma}}\right)\frac{\pi}{2}$.
Отже, при $\alpha_{0}\geq\frac{2}{\sqrt{\gamma}}$

$$I_{\gamma}\leq p^{\gamma-1}\cdot\frac{9}{4^{\frac{8}{3}}}\left(|a_{1}-a_{2}|\cdot|a_{1}-a_{3}|\cdot|a_{2}-a_{3}|\right)^{\frac{2}{3}}\leq$$
$$\leq4p_{0}^{\gamma-1}\cdot\frac{9}{4^{\frac{8}{3}}}\sin^{\frac{2}{3}}\frac{\pi}{\sqrt{\gamma}}\cdot\sin^{\frac{4}{3}}\left(1-\frac{1}{\sqrt{\gamma}}\right)\frac{\pi}{2}.$$

Підставивши відповідні значення $\gamma$ і $p_{0}$, отримаємо, що
$I_{\gamma}\leq0.2936< I_{1.7}^{0}$, тобто при
$\alpha_{0}\geq\frac{2}{\sqrt{\gamma}}$ конфігурації областей не
можуть бути екстремальними.

Звідси, $\alpha_{0}\leq\frac{2}{\sqrt{\gamma}}$, і ми можемо
застосовувати нерівність (2).

Далі, використовуючи результат, отриманий під час доведення теореми
5.2.3 із [9], запишемо наступну нерівність:

\begin{center}
$I_{\gamma}\leq\frac{1}{\sqrt{\gamma}} \cdot[\prod\limits_{k=1}^2
2^{\sigma_{k}+6}\cdot\sigma_{k}^{\sigma_{k}+2}\cdot(2-\sigma_{k})^{-\frac{1}{2}(2-\sigma_{k})^{2}}(2+\sigma_{k})^{-\frac{1}{2}(2+\sigma_{k})^{2}}]^{\frac{1}{2}}$
\end{center}

де $\sigma_{k}=\sqrt{\gamma}\cdot\alpha_{k}$ Введемо функцію

\begin{center}
$\Psi(\sigma)=
2^{\sigma+6}\cdot\sigma^{\sigma+2}\cdot(2-\sigma)^{-\frac{1}{2}(2-\sigma)^{2}}(2+\sigma)^{-\frac{1}{2}(2+\sigma)^{2}}$
\end{center}

для $\sigma\in[0,2]$ і, використовуючи її поведінку на цьому
проміжку, доведемо екстремальність конфігурації областей $D_{0},
D_{1}, D_{2}$.

$\Psi(\sigma)$ логарифмічно опукла на проміжку $[0;x_{0}]$, де
$x_{0}\approx1.32$. На проміжку $[0;x_{1}]$, $x_{1}\approx1.05$
функція зростає від $\Psi(0)=0$ до $\Psi(x_{1})\approx0.9115$,
спадає на проміжку $[x_{1};x_{2}]$ $x_{2}\approx1.6049$до
$\Psi(x_{2})\approx0.86$, а на проміжку $[x_{2};2]$ зростає до
$\Psi(2)=1$. Точка, $x_{3}$ для якої
$\Psi(x_{3})=\Psi(x_{1})\approx0.9115$ $x_{3}\approx 1.9$. Тепер,
використовуючи рівність $\sigma_{1}+\sigma_{2}=2\sqrt{\gamma}$,
доведемо, що
$\Psi(\sigma_{1})\cdot\Psi(\sigma_{2})\leq(\Psi(\sqrt{\gamma}))^{2}\approx0.8308$.
Для $\sigma_{k}\in[0;x_{0}]$ відповідний висновок робимо із
логарифмічної опуклості функції $\Psi(\sigma)$. Для
$\sigma_{2}\in[x_{0};x_{3}]$ із властивостей графіка функції
$\Psi(\sigma)$,  $\Psi(\sigma_{2})\leq\Psi(\sqrt{\gamma})$ і
$\Psi(\sigma_{1})\leq\Psi(\sqrt{\gamma})$, а тому
$\Psi(\sigma_{1})\cdot\Psi(\sigma_{2})\leq(\Psi(\sqrt{\gamma}))^{2}$.
Якщо ж $\sigma_{2}\in[x_{3};2]$, то $\Psi(\sigma_{2})<\Psi(2)=1$,
$\Psi(\sigma_{1})<\Psi(0.2)\ll0.4$, то звідси
$\Psi(\sigma_{1})\cdot\Psi(\sigma_{2})<0.4<(\Psi(\sqrt{\gamma}))^{2}$.
Отже, $I_{\gamma}\leq I_{\gamma}^{0}(x_{1})$, тому екстремальних
конфігурацій областей ми не отримаємо.

Для $\gamma=1.7$ теорему доведено. Доведемо, що функціонал
$I_{\gamma}^{0}$ як функція від $\gamma$ монотонно спадає на
проміжку $[1.5;1.7]$. Для цього візьмемо логарифмічну похідну від
виразу

\begin{center}
$I_{\gamma}^{0}=\left(\frac{4}{n}\right)^{n}\cdot\frac{\left(\frac{4\gamma}{n^{2}}\right)^{\frac{\gamma}{n}}}{\left(1-\frac{\gamma}{n^{2}}\right)^{n+\frac{\gamma}{n}}}\cdot\left(\frac{1-\frac{\sqrt{\gamma}}{n}}{1+\frac{\sqrt{\gamma}}{n}}\right)^{2\sqrt{\gamma}}$
\end{center}

при $n=3$. Отримаємо

\begin{center}
$(\ln(I_{\gamma}^{0}))'=\frac{1}{3}\ln\frac{4\gamma}{9}-\frac{1}{3}-\frac{\gamma}{3}\cdot\ln(1-\frac{\gamma}{9})+\frac{9+\gamma}{27-3\gamma}+\frac{1}{\sqrt{\gamma}}\ln(1-\frac{\sqrt{\gamma}}{3})-\frac{1}{3-\sqrt{\gamma}}-\frac{1}{\sqrt{\gamma}}\ln(1+\frac{\sqrt{\gamma}}{3})-\frac{1}{3+\sqrt{\gamma}}$
\end{center}

На проміжку $[1.5;1.7]$

$$(\ln(I_{\gamma}^{0}))'\leq\frac{1}{3}\ln\frac{6.8}{9}-\frac{1}{3}-\frac{1.7}{3}\ln(1-\frac{1.7}{9})+\frac{9+1.7}{27-3\cdot1.7}-\frac{1}{3-\sqrt{1.5}}\approx$$
$$\approx-0.0934-\frac{1}{3}+0.1186+0.4886-0.5633\approx-0.3828<0$$ тому дана функція є монотонно спадною, а отже
$I_{\gamma}^{0}>I_{1.7}^{0}$ при $\gamma\in[1;1.7]$. За
властивостями функції $\sin x$ отримаємо

\begin{center}
$p^{\gamma-1}\cdot\frac{9}{4^{\frac{8}{3}}}\left(|a_{1}-a_{2}|\cdot|a_{1}-a_{3}|\cdot|a_{2}-a_{3}|\right)^{\frac{2}{3}}$
\end{center}

Таким чином
$\frac{I_{\gamma}}{I_{\gamma}^{0}}\leq\frac{I_{1.7}}{I_{1.7}^{0}}<1$.
Звідси для $\gamma\in(1.5;1.7]$ $I_{\gamma}\leq I_{\gamma}^{0}$, а
тому $I_{\gamma}^{0}$ - шукана екстремальна конфігурація областей.

Теорему доведено.

 \begin{center}

{\bf СПИСОК ВИКОРИСТАНИХ ДЖЕРЕЛ}

\end{center}

\end{document}